\documentclass{amsart}

\newtheorem{theorem}{Theorem}[section]
\newtheorem{proposition}[theorem]{Proposition}
\newtheorem{lemma}[theorem]{Lemma}
\newtheorem{corollary}[theorem]{Corollary}
\newtheorem{conjecture}[theorem]{Conjecture}

\theoremstyle{definition}
\newtheorem{definition}[theorem]{Definition}

\theoremstyle{remark}

\numberwithin{equation}{section}

\newcommand{\inv}{^{-1}}
\newcommand{\del}{\partial}
\newcommand{\delb}{\overline{\partial}}

\newcommand{\bbC}{{\mathbb C}}
\newcommand{\bbH}{{\mathbb H}}
\newcommand{\bbZ}{{\mathbb Z}}
\newcommand{\bbR}{{\mathbb R}}
\newcommand{\bbP}{{\mathbb P}}

\newcommand{\calW}{{\mathcal W}}
\newcommand{\calS}{{\mathcal S}}
\newcommand{\calM}{{\mathcal M}}
\newcommand{\calN}{{\mathcal N}}
\newcommand{\calB}{{\mathcal B}}

\newcommand{\calH}{{\mathcal H}}
\newcommand{\calE}{{\mathcal E}}
\newcommand{\calF}{{\mathcal F}}

\newcommand{\calV}{{\mathcal V}}
\newcommand{\calG}{{\mathcal G}}

\DeclareMathOperator{\tr}{Tr}

\newcommand{\dbar}{\bar\partial}
\newcommand{\norm}[1]{\left\Vert #1 \right\Vert}

\newcommand{\dtheta}{\partial_\theta}

\begin{document}

\title{Nearly K\"ahlerian Embeddings of Symplectic Manifolds}
\author{David Borthwick}
\address{Department of Mathematics and Computer Science\\
Emory University\\Atlanta, GA 30322}
\email{davidb@mathcs.emory.edu}
\thanks{First author supported in part by NSF grant DMS-9796195 
and by an NSF postdoctoral fellowship.}
\author{Alejandro Uribe}
\address{Mathematics Department\\
University of Michigan\\Ann Arbor, MI 48109-1109}
\email{uribe@math.lsa.umich.edu}
\thanks{Second author supported in part by NSF grant DMS-9623054.}
\date{December 4, 1998}

\maketitle

\tableofcontents

\newcommand{\Lk}{L^{\otimes k}}

\section{Introduction and statement of results}

Let $X,\omega$ be a $2n$-dimensional
compact symplectic manifold without boundary.
There always exists a compatible almost-complex structure $J$ on $X$,
i.e. 
$$
\omega(Ju,Jv) = \omega(u,v),\hskip.5in \omega(\cdot,J\cdot) \gg 0.
$$
Let $g$ denote the associated Riemannian metric on $X$, that is
$g(u,v)\,=\,\omega(u,Jv)$.   Assume that $\omega$ defines an integral
cohomology class, which means that there exists a complex
line bundle, $L\to X$, with a connection $\nabla$ whose curvature
is $-i\omega$.\footnote[1]{This convention is fairly common in complex geometry.  
But it is also common to include a factor of $2\pi$, and in quantum mechanics there
is yet a different convention.  The choice affects our normalization of the 
Fubini-Study form, and will be responsible for certain factors of 2 in Sections
4 and 6.}  In \cite{BU} we showed how
to ``quantize" $(X,\omega, J)$ by a family of finite-dimensional
Hilbert spaces, 
$\calH_k\subset C^\infty(X,\Lk)$, defined for all $k$ sufficiently large,
and we observed that the usual procedure of algebraic geometry yields
well-defined maps,
\begin{equation}\label{uno}
F_k :\, X\longrightarrow \bbP\calH_k^*.
\end{equation}
The $\calH_k$ are the span of certain eigensections of the Laplacian on 
$C^\infty(X,\Lk)$ defined by $g$ and $\nabla$;  see \S2 for the precise
definition.  The dimension of $\calH_k$ is the Riemann-Roch polynomial of
$X$ evaluated at $k$, that is 
\[
\mbox{dim}\ \calH_k\,=\,
\int\,e^{k[\omega ]}\,\tau\,=\, \mbox{Vol}(X)k^n + O(k^{n-1}),
\]
where $\tau$ is the Todd class. 

In this paper we study the asymptotic geometry of the
maps $F_k$.  We propose the general philosophy that
the quantization $\{\calH_k\}$ is a natural object that can be used
to study the geometry of $(X,\omega)$.  The main results are as follows:

\begin{theorem}\label{Main}
Let $\Omega_k$ and $g^0_k$ denote respectively the symplectic form and
K\"ahler metric of $\bbP\calH_k^*$.  Then there exist positive constants
$C_1$ and $C_2$ and $k_0$ such that for all $k>k_0$
\begin{equation}\label{dos}
\forall x\in X\qquad |\frac{1}{k}F_{k}^{*}\Omega_{k}- \omega|_{x}\leq 
\frac{C_{1}}{k}
\end{equation}
and
\begin{equation}\label{tres}
\forall x\in X\qquad |\frac{1}{k}F_{k}^{*}g^{0}_{k}-g|_{x}\leq 
\frac{C_{2}}{k}\,.
\end{equation}
Moreover, $F_k$ is injective for large $k$ and therefore an embedding.
\end{theorem}
The precise meaning e.g. of (\ref{tres}) is: $\forall u\in T_xX$
\begin{equation}\label{maya}
\norm{ dF_k(u)}^2_{g_k^0}\,=\,k\,\norm{u}^2_g + O(1),
\end{equation}
where the estimate is uniform on the unit ball bundle, $\{ \norm{u}\leq 1\}$.
The following is an easy consequence:

\begin{corollary}\label{Co}
The $F_k$ are approximately pseudo-holomorphic
for large $k$.  More precisely:
\begin{equation}\label{don}
\frac{1}{k} \|\partial(F_k)_x\|\,=\,1+ O(k^{-1}),\quad\mbox{while}\quad
\frac{1}{k}\|\overline{\partial}(F_k)_x\|\,=\ O(k^{-1}),
\end{equation}
where $\|\ \|$ denotes the operator norm, and the estimates are uniform in $x\in X$.
\end{corollary}

Theorem \ref{Main} says that the $F_k$ are projective embeddings which,
for $k$ large, are both `nearly' symplectic and `nearly' isometric. 
(Recall that by a Theorem of Tischler, \cite{T},
closed integral symplectic manifolds admit symplectic embeddings
into projective space, but in his proof
one has no information about the isometry properties of the embedding.)
In case $X$ is K\"{a}hler estimates as above were obtained by Tian, \cite{Ti}.
Recently, Zelditch has re-proved Tian's results using microlocal methods,
\cite{Ze}.  Our proof of Theorem (\ref{Main}) is in part parallel to Zelditch's.
However we have to work harder since we don't have the general 
benefits of working in the holomorphic category (e.g.
$F_k$ is not automatically holomorphic).

The proof of Theorem \ref{Main} relies heavily on the machinery of
Fourier integral operators of Hermite type, of Boutet de Monvel and
Guillemin \cite{BG}.  The strategy is to realize the family of spaces
$\calH_k$ as the Fourier decomposition  of the image of a generalized
Szeg\"o projector $\Pi$. Theorem \ref{Main} is then essentially a
statement about the symbol of $\Pi$. Our proof will be purely symbolic,
so we could replace $\Pi$ by any generalized Szeg\"o projector with the
same microlocal properties.  
 
It should also be
possible to twist $L$ with a fixed Hermitian vector
bundle, $\calE\to X$, and work with generalized Szego
projectors acting on the sections of a vector bundle.  
Just as in algebraic geometry, one then obtains embeddings into Grassmannians.

One of our goals in studying the maps $F_k$ is to ultimately give a
microlocal  proof of the recent ground-breaking theorem of Donaldson,
\cite{D}, on the  existence of symplectic submanifolds of a
symplectic manifold that are Poincar\'e dual to $k[\omega]$ for $k$ large.
\begin{conjecture}
For large $k$ the transverse hyperplane sections of
$F_k$ are symplectic submanifolds of $X$.
\end{conjecture}
This does not follow immediately from Theorem \ref{Main}.  We need to know
in addition
that for large $k$ there are hyperplane sections which intersect $F_k(X)$ 
with a minimum angle $\epsilon>0$ which is independent of $k$.  This
is a quantitative transversality issue, just as in Donaldson's proof.
His arguments however do not apply directly here: his model sections of $\Lk$
can be said to be close to our space $\calH_k$ only in an $L^2$ sense.  
In any case, we hope that the microlocal estimates on $F_k$ will eventually yield
an independent proof.

Section 2 is devoted to preliminaries, and in \S 3 we review some
facts about Hermite FIOs.
Theorem \ref{Main} is proved in \S 4, with some details relegated to
an Appendix.   In \S5 we give an example of the kind of weak version
of Donaldson's theorem which does follow easily from Theorem \ref{Main}.
Finally, in \S6 we consider the relationship between classical dynamics
on $X$ and quantum dynamics on $\bbP\calH_k^*$.

\section{Preliminaries}

\subsection{Almost-K\"ahler quantization.}

Henceforth $X,\,\omega,\, J,\, L,\,\nabla$ will be as in \S 1.
The Hermitian metric and connection on $L$ induce corresponding
structures on $\Lk$.  Together with the (fixed) metric on $X$,
this defines a Laplace operator, $\Delta_{k}$, acting on
$C^{\infty}(X,\Lk)$.  We now define the {\em rescaled Laplacian},
\begin{equation}\label{rlap}
\calB_{k}\, :=\,\Delta_{k}- n\,k.
\end{equation}
In \cite{GU} it was observed that Mellin's inequality implies the 
existence of constants $C_{1},\ C_{2} >0$ such that the spectrum
of $\calB_{k}$ is contained in
\[
(-C_1\,,\, C_1)\ \bigcup\ (k\,C_{2}\,,\infty).
\]
Therefore for large $k$ the spectrum exhibits a gap of size $O(k)$.
Moreover, for large $k$ the number of eigenvalues (with multiplicities)
of $\calB_{k}$ in $(-C_1\,,\, C_1)$ is precisely equal to the
value of the Riemann-Roch polynomial at $k$. 
We can now define the spaces appearing in Theorem (\ref{Main}):

\begin{definition}
For all sufficiently large $k$, the quantizing space
$\calH_{k}\subset C^{\infty}(X,\Lk)$ is defined to be the span of those 
eigensections of $\calB_{k}$ with eigenvalues in $(-C_1\,,\, C_1)$.
\end{definition}

By the Bochner-Kodaira formula (e.g. \cite{BGV} Proposition 3.71),
in case $(X,\omega,J)$ is a K\"ahler manifold and $L$ holomorphic 
Hermitian with the induced connection, the operator $\calB_{k}$
is precisely the $\overline{\partial}$-Laplacian and $\calH_{k}$ is
the space of the holomorphic sections of $\Lk$.

\subsection{The unit circle bundle.}
Let $Z\subset L^{*}$ denote the unit circle bundle.  The connection
on $L$ induces one on $Z$ (as a principal $S^{1}$ bundle), and this
together with the Riemannian metric on $X$ and the standard
metric on $S^{1}=\bbR/2\pi\bbZ$ induces a metric on $Z$ such that
the projection $Z\to X$ is a Riemannian submersion with totally 
geodesic fibers.  The Laplacian on $Z$, $\Delta_{Z}$, commutes with
$\partial_{\theta}$, the infinitesimal generator of the circle 
action, and therefore it also commutes with the ``horizontal 
Laplacian'',
\begin{equation}
\Delta_{h}\,:=\, \Delta_{Z} - \partial_{\theta}^{2}.
\end{equation}
As is well-known, the decomposition of $L^{2}(Z)$ into $S^{1}$ 
isotypes, $\displaystyle{L^{2}Z = \oplus_{k}L^{2}(Z)_{k}}$ is such 
that for each $k$ $L^{2}(Z)_{k}$ is naturally isomorphic with
$L^{2}(X,L^{k})$.  Under this isomorphism $\Delta_{h}$ gets identified
with the Laplacian on sections of $L^{k}$ (induced by the metrics and the
connection), and therefore the rescaled
Laplacian of the previous subsection is induced by the single 
operator on $Z$,
\begin{equation}
\calB\,:=\,\Delta_{h} - n\partial_{\theta}\,.
\end{equation}
It will be convenient to identify $\calH_{k}$ with a subspace of
$L^{2}(Z)_{k}$, and let 
\begin{equation}\label{da}
\calH\,:=\,\widehat{\bigoplus_{k}}\,\calH_k \qquad \mbox{and}\qquad
\Pi: L^2(Z)\to\calH
\end{equation}
the orthogonal projection.
In Appendix B we show that $\Pi$
is a Fourier integral operator of Hermite type.  

\subsection{The coherent state map.}
Let $\Pi_{k}(x,y)$ denote the Schwartz kernel of
the orthogonal projection,
\begin{equation}
\Pi_{k}:\, L^{2}(Z)\to \calH_{k}\,.
\end{equation}
\newcommand{\Psik}{\Psi_{k}}

\begin{definition} (\cite{Raw})
The \textit{coherent state} in $\calH_{k}$ associated to a point $p\in Z$
is 
\[
\Pi_{k}(\cdot,p) \in \calH_k.
\]
And the coherent state map is 
\begin{equation}\label{a0}
\begin{array}{rcc}
\Psik :\,Z&\to &\calH_{k}\\
p & \mapsto & \Pi_k(\cdot,p)
\end{array}
\end{equation}
\end{definition}

\noindent
{\bf Remarks:}
\par\noindent
{\bf 1.} By their definition the coherent states have the \textit{reproducing property}
\begin{equation}\label{a1}
\forall f\in\calH_k,\ p\in Z \qquad f(p)\,=\,\langle\Psik(p)\,,\,f\rangle.
\end{equation}

\medskip\noindent{\bf 2.}
The map (\ref{a0}) clearly induces the map (that we'll denote by $\Psik^o$) 
\begin{equation}\label{a3}
\begin{array}{rcl}
\Psik^o: X & \to & \bbP\calH_k\\
x &\mapsto &  [\Psik(p_x)]
\end{array}
\end{equation}
where $p_{x}\in Z$ is any point projecting to $x$ and
$[\Psik(p)]$ denotes the complex line through $\Psik(p)$.

\medskip\noindent{\bf 3.}
Recall that the standard map (\ref{uno}) is defined by:
\[
F_k(x)\,=\,\{\,f\in\calH_k\ :\ f(p_x) = 0\ \mbox{where}\ p_x\in P\ 
\mbox{satisfies}\ \pi(p_x)=x\,\},
\]
where one identifies an element of $\bbP\calH^*$ with a hyperplane
in $\calH$.   Let 
\[
j:\ \bbP\calH_k \to \bbP\calH^*_k
\]
be the (anti-holomorphic) map defined by the Hermitian metric on $\calH_k$.
Then by the reproducing property (\ref{a1}) the diagram
\begin{equation}\label{a4}
\begin{array}{rcc}
 & & \bbP\calH^*_k \\
 & \stackrel{F_k}{\nearrow} & \uparrow\ j\\
X & \stackrel{\Psik^o}{\longrightarrow} & \bbP\calH_k
\end{array}
\end{equation}
commutes.  Since the statements of Theorem \ref{Main} are invariant
under composition of $F_k$ by $j\inv$, it suffices to prove the Theorem
for the map $\Psik^o$.

\section{Hermite FIOs.}

\subsection{Structure of $\Pi$}
The key technical feature of the almost-K\"ahler quantization 
is the fact that the orthogonal projector 
$\Pi:L^2(Z)\to\calH$ is an FIO of Hermite type.  This means
that locally we can represent its kernel as an oscillatory integral
of a particular form.  

Given a manifold $M$ and a conic isotropic submanifold
$\Sigma\subset T^{*}M\setminus\{ 0\}$,  we'll denote by $J^{m}(M,\Sigma)$
the spaces of Hermite distributions introduced by Boutet de Monvel and
Guillemin in \cite{BG}.  Hermite FIOs have kernels in this class, which
we'll describe in detail below.
The main technique of \cite{GU} was to establish that $\calB$ can
be modified by an pseudodifferential operator $R$ of order zero so that
the projector onto $\ker(\calB-R)$ is such an Hermite FIO.  In fact, this
implies that our projector $\Pi$ is also Hermite:
\begin{theorem}\label{Bg}
Let
\begin{equation}\label{1d}
\Sigma\,=\,\{\,(p,\,r\alpha_{p};
p,\,-r\alpha_{p})\ ;\ r>0,\ p\in Z\,\} \subset T^{*}(Z\times
Z)\setminus \{ 0\},
\end{equation}
where $\alpha$ is the connection form on $Z$. 
Then $\Pi\in J^{1/2}(Z\times Z, \Sigma)$. 
\end{theorem}

\noindent
As the proof is somewhat technical, we defer it to Appendix A.

A distribution $u\in J^m(Z\times Z, \Sigma)$ has the following
local description. We can cover $Z\times Z$ by (finitely many) coordinate
patches, so that in patch the kernel is written as an integral over phase
variables
$(\tau,\eta)\in\bbR_+\times\bbR^{2n}\backslash\{0\}$:
\begin{equation}\label{ulocal}
u(x,y) = \int e^{i\tau f(x,y) + i\eta\cdot g(x,y)}\>
a(x,y,\tau,\eta/\sqrt{\tau})\>d\tau\>d\eta.
\end{equation}
Here the phase functions $f$ and $g_j$, $j=1,\dots,2n$, are
smooth functions which ``parametrize'' $\Sigma$ in the sense that
$$
\Sigma = \{(x,y;\>\tau df(x,y)): \tau>0,\; f(x,y) = g_1(x,y) = \dots =
g_{2n}(x,y) = 0\}.
$$
Thus $f=g_1=\dots=g_{2n}=0$ must define the diagonal in $Z\times Z$,
and $df(x,x) = (\alpha_x,-\alpha_x) \in T^*(Z\times Z)$.
The amplitude, $a(x,y,\tau,\xi)$, is polyhomogeneous in $\tau$ with decreasing
half-integer powers, and rapidly decreasing as a function of $\xi$.
For a distribution in $J^m(Z\times Z;\Sigma)$, the leading term 
in the amplitude will be of the form $\tau^{m-1/2} a_0(x,y,\xi)$.
An application of stationary phase shows that these requirements insure
that the wave front set of (\ref{ulocal}) is contained in $\Sigma$.

The symbol of a distribution in $J^m(Z\times Z;\Sigma)$ is the invariant
object corresponding to $a_0(x,y,\xi)$.  This object has two parts.  (As
usual, to get an invariant symbol one should consider distributions acting on
half-forms.)  The first part of the symbol is a half-form on
$\Sigma$, i.e. the square-root of a top-degree differential form. (This
would be the full symbol for a Lagrangian distribution.) 
To make the square-root globally well-defined, we need a metalinear
structure on $\Sigma$.  

The second part of the symbol captures the $\xi$ dependence of $a_0$.  
Let $\Sigma\subset T^{*}M\setminus\{ 0\}$ be an isotropic submanifold.
At each point $z\in \Sigma$ we define the symplectic normal space:
\[
\calN_z\,:=\,T_z\Sigma^\bot/T_z\Sigma,
\]
and $\calN=\bigcup_{z\in\Sigma}\calN_z$ the symplectic normal bundle
of $\Sigma$. 
Suppose we are given a metaplectic structure on $\calN_z$. Then for every
$z\in \Sigma$ one can construct the Heisenberg group
$\bbH_z = \calN_z\oplus\bbR$ associated with $\calN_z$ and the metaplectic 
group $\calM_z$, the double cover of the group of linear symplectic
transformations of $\calN_z$.  Let $\calS_z$ be the space of smooth
vectors in the metaplectic representation of $\calM_z$.  The local
representative $a_0(x,y,\xi)$ corresponds to an element of $\calS_z$.

The resulting symbol is thus a smooth half-form on 
$\Sigma$ with values in the bundle
$\calS:=\bigcup_{z\in\Sigma}\calS_z$, called a {\it symplectic spinor}.
The symbol will be homogeneous with respect to the $\bbR_+$
action on $\Sigma$, reflecting the $\tau$ dependence of the leading
amplitude.

The structures necessary to define the symbol map always exists in our
case.  In particular, the manifold
$Z$ always possesses a metalinear structure (Lemma 2.8 of \cite{BPU1}).
From this choice we can derive the product of a metalinear 
structure on $\Sigma$ and a metaplectic structure on $\calN$.  

\begin{lemma}
The symplectic normal space to (\ref{1d}) at $z=(p,\alpha_p;
p,-\alpha_p)$ is naturally isomorphic with 
\begin{equation}\label{iso}
\calN_z\,\cong \,H_p\oplus H_p
\end{equation}
where $H_p\subset T_pZ$ is the horizontal space at $p$.
Therefore 
\begin{equation}\label{db}
\calS_z\,\cong \,\calS_p\otimes\calS_p
\end{equation}
where $\calS_p$ denotes the space of smooth vectors of the
metaplectic representation of the metaplectic group of $H_p$. 
\end{lemma}
\begin{proof}
The isomorphism (\ref{iso}) is as follows.  First it is clear that
\[
\calN_z\,\cong\, E_p\oplus E_p,
\]
where $E_p$ is the symplectic orthogonal  to the tangent space at
$(p,\alpha_p)$ to the {\em symplectic} submanifold
\[
\{\,(p,r\alpha_p)\;;\; p\in Z,\ r>0\,\}\subset T^*Z.
\]
Then one can check that the differential of the natural projection
$T^*Z\to Z$ induces a symplectomorphism $E_p\to H_p$ 
\end{proof}

Technically, in order to define the symbol of $\Pi$ invariantly
we should consider $\Pi$ as an operator on half-forms.
We can circumvent this issue, however, because
if $\Sigma$ is as in (\ref{1d}) then one has a natural identification
\[
\Sigma\,\cong\, Z\times\bbR^{+}.
\]
Under the identification $(p,r)\mapsto (p,r\alpha_p) \in T^*Z$,
this space is symplectic.
Therefore $\Sigma$ possesses a natural nowhere-vanishing
half-form.  Using this half-form and (\ref{db}) we can
conclude:

\begin{lemma}\label{TeS}
If $u$ is an Hermite distribution on $Z\times Z$ associated to
(\ref{1d}), its symbol at $z=(p,\alpha_p ; p,-\alpha_p)$ can be 
naturally identified with an operator
\[
\sigma_u(p)\, : \calS_p\to \calS_p.
\]
\end{lemma}

For case of the Szeg\"o projector (the K\"ahler case), the symbol of $\Pi$
was calculated in \cite{BG}.  It is determined by the K\"ahler
structure on the horizontal fibers $H_p\subset T_{p}Z$.  The
integrability of the almost complex structure plays no role in the symbol, 
and we get the same result for the almost-K\"ahler case (proven in Appendix A):

\begin{proposition}\label{szsym}
The symbol of $\Pi$ at $(p,\,\alpha_{p}; p,\,-\alpha_{p})$
is the rank-one projector
$|e><e|$, where $e$ is the ground state of the
harmonic oscillator associated to the K\"ahler structure of the
horizontal subspace $H_p\subset T_{p}Z$.
\end{proposition}

\subsection{Asymptotics}

Given $u\in J^m(Z\times Z, \Sigma)$, we can decompose 
into isotypes:  $\displaystyle{ u\,=\,\oplus_k\,u_k }$, where
$$
u_k(x\cdot e^{i\theta}, y)\,=\, e^{ik\theta}\,u_k(x,y).
$$
The wave front set $\Sigma$ consists of a ray in the cotangent
bundle over each point in the diagonal.  This ray should be thought
of as the ``large $k$'' direction.  For if $J$ is the Hamiltonian
generating the $S^1$ action on $T^*Z$, we have $J(x,r\alpha_x) = r$.
A simple wave front set argument shows that $u_k$ is a smooth
function for each $k$.  So knowledge of the singularities of $u$ may be
translated into information on the large $k$ behavior of $u_k$.

To make this statement precise, we use the method of stationary phase.
\begin{lemma}\label{spest}
For $u\in J^m(Z\times Z, \Sigma)$, decomposed as above,
there exists a constant $C$ such that
$$
\sup_{Z\times Z} |u_k| \le Ck^{m+n-1/2}.
$$
\end{lemma}
\begin{proof}
To pick off the $k$-th isotype of $u$ we average over
$\theta$:
$$
u_k(x,y) = \int_0^{2\pi} e^{-ik\theta}u(x\cdot e^{i\theta}, y) \>d\theta.
$$
Since $Z\times Z$ is compact, it suffices to consider the local form of
$u$ given by (\ref{ulocal}).  We will focus on the $\tau$ and $\theta$
integrations for the moment.  So for the leading term we study
$$
w_k(x,y,\eta) = \int e^{-ik\theta} e^{i\tau f(x\cdot e^{i\theta}, y)}
\tau^{m-1/2} a_0(x\cdot e^{i\theta}, y,\eta/\sqrt{\tau}) \>d\tau\> d\theta. 
$$
Since this expression contains local cutoffs, we can extend the $\theta$
integration to $\bbR$ to simplify the notation.  
Rescaling the integral by $\tau\to k\tau$ yields
\begin{equation}\label{spresult}
w_k(x,y,\eta) = k^{m+1/2}
\int e^{-ik\theta + ik\tau f(x\cdot e^{i\theta}, y)}
\tau^{m-1/2} a_0(x\cdot e^{i\theta}, y,\eta/\sqrt{k}) \>d\tau\> d\theta. 
\end{equation}

The phase function $f$ must satisfy the condition $df|_{(x,x)} =
(\alpha_x, -\alpha_x)$.  Thus $\frac{d}{d\theta} f(x\cdot e^{i\theta}, x)
= \alpha_x(\del_\theta) = 1$.  There is freedom in the choice of $f$,
so that we can require that this property hold off the diagonal as well.
In fact we may as well choose $f$ so that $f(x\cdot e^{i\theta}, y) = \theta +
f(x,y)$.

With this assumption, we are ready to apply stationary phase to the
$\tau,\theta$ integration.  The stationary
point in the phase $\tau f(x\cdot e^{i\theta}, y) - \theta$ occurs at
$(\tau, \theta) = (1,- f(x,y))$.   Applying  a stationary phase estimate with
parameters, such as Theorem 7.7.6 of
\cite{H}, we obtain a uniform estimate
$$
w_k(x,y,\eta) = C k^{m-1/2} e^{ikf(x,y)} a(x,y,\eta/\sqrt{k\tau})
+ O(k^{m-3/2}).
$$
Finally, we plug this back into the $\eta$ integral and rescale to remove
the $k$, which adds a power $k^n$.  The lower order terms in the expansion of
the amplitude are handled similarly.   
\end{proof}

In the above estimate, by ignoring the $\eta$ integration we did not
distinguish between the behavior on and off the diagonal.  If include the
$\eta$ integration and restrict
to a compact subset of $P\times P$ whose projection to $X\times X$ doesn't
intersect the diagonal, then the phase has no stationary points of the phase, 
and we have a uniform estimate $u_k(x,y) = O(k^{-\infty})$.

We also will need to understand precisely the leading behavior of
$u_k$ on the diagonal.  This could be established by an
argument as above, but instead we will appeal to more general results
from \cite{BG}.   In this argument the stationary phase approximation
is implicit in the composition theorem cited.

\begin{lemma}\label{Tech}
For $u\in J^m(Z\times Z, \Sigma)$, decomposed as above,
one has an asymptotic expansion
\begin{equation}\label{tech}
u_k(p,p)\,\sim\, \frac{k^{n+m-\frac{1}{2}}}{(2\pi)^n}\,
\sum_{j=0}^\infty\, k^{-j/2}\,f_j(p),
\end{equation}
as $k\to\infty$.  The coefficients $f_j$ are smooth and the expansion
is in the $C^\infty$ topology.
Moreover, using the notation of Lemma \ref{TeS}, the leading coefficient 
is 
\begin{equation}\label{tech2}
f_0(p)\,=\,\tr\sigma_u(p),
\end{equation}
where $\tr\sigma_u(p)$ denotes the integral along the
diagonal of the Schwartz kernel of $\sigma_u (p)$.

If $u$ in addition has definite parity, as defined in \S9 of \cite{BG},
the asymptotic expansion goes down by integral powers of $k$.
\end{lemma}
\begin{proof}
The following technique has been used several times before, 
\cite{BPU1}, \cite{BPU2}.
Let $P: C^\infty(Z\times Z) \to C^\infty(S^1)$ be given by 
\begin{equation}\label{pfio}
Pu(\theta) =
u(p\cdot e^{i\theta}, p),
\end{equation}
where $p$ is fixed for now, so that
$$
Pu(\theta) = \sum_k e^{ik\theta} u_k(p,p).
$$
Clearly $P$ is an ordinary Fourier integral operator.
One can check that the main composition theorem of \cite{BG}, Theorem 7.5,
applies so that $P(u)$ is a distribution on $S^1$ conormal to the
identity. Moreover,
Theorem 9.9 of \cite{BG} applies to show that $P$ preserves parity.  (The
excess of the composition fiber product is $2n$, and $\dim (Z\times Z) - 2n$
is even.)  This implies the existence of the asymptotic expansions of $u_k(p,p)$
with the desired properties, since these are the Fourier coefficients of $P(u)$.  
The leading coefficient in the expansion of $u_k(p,p)$ is
given by the symbol of $P(u)$.  We defer this symbolic calculation
to Appendix B.

The same argument goes through if we consider $P$ as an operator
$P: C^\infty(Z\times Z) \to C^\infty(S^1\times Z)$,
so the asymptotic expansion is in the $C^\infty$ topology. 
\end{proof}

\section{Nearly K\"ahlerian embedding}

In this section we'll prove Theorem \ref{Main}.
As noted in \S2, the statements about $F_k$ are equivalent
to those for $\Psik^o$.  The microlocal structure described in the
preceding section gives us direct estimates of the coherent state
map $\Psik$, so we will work mainly with this.

We must compute the the pullback by $F_k$ of the Fubini-Study
hermitian metric on $\bbP\calH_k$.  We'll denote the pullback at the point
$x\in X$ by $h_{x}(\cdot,\cdot)$, i.e.
$$
h_x = F_k^*g_k^0 + iF_k^*\Omega_k
$$
To compute this, we lift tangent
vectors horizontally to $T_pP$, where $\pi(p) = x$, pushforward by
$\Psik$, and then plug the resulting vectors into the homogeneous coordinate version
of the Fubini-Study hermitian metric.  Letting $u_h, v_h$ denote the horizontal lifts,
\begin{equation}\label{hxdef}
\frac{1}{2}h_x(u,v) = \frac{\langle d(\Psik)_p(v_h),
d(\Psik)_p(u_h)\rangle}{\norm{\Psik(p)}^2}
- \frac{\langle d(\Psik)_p(v_h), \Psik(p)\rangle  
\langle \Psik(p), d(\Psik)_p(u_h)\rangle}{\norm{\Psik(p)}^4}.
\end{equation}
Here the inner products and norms are those of $\calH_k$ (or equivalently
$L^2(Z)$), and the normalizing factor of $2$ is consistent with our choice
that the curvature of the quantizing line bundle is $(-i)$ times the symplectic
form.

\subsection{Calculation of the pullback.}

Let $\nu_k$ be the restriction of $\Pi_k$ to the diagonal,
$$
\nu_k(p) := \Pi_k(p,p).
$$
This is in fact the same as the norm square of the coherent states, 
by the reproducing property:
$$
\nu_k(p)\,=\,\norm{\Psik(p)}^2.
$$
(Incidentally, this proves that $\nu_k$ is strictly positive.)
In the K\"ahler case, the function $\nu$ was introduced by
Rawnsley in \cite{Raw} (and further studied in \cite{CGR}, denoted by
$\theta$ there), where it
is shown that its constancy is equivalent to the coherent state map
being symplectic.  This doesn't hold in general,
not even for all K\"ahler manifolds.
However, we can show that $\nu_k$ is at least asymptotically constant, even
in the almost-K\"ahler case.

To proceed we calculate and estimate the terms appearing in $h_x$.
\renewcommand{\u}{u}
\renewcommand{\v}{v}
\newcommand{\kpi}{\Bigl(\frac{k}{2\pi}\Bigr)}
\begin{theorem}\label{Lapapa}
For $p\in Z$, there is an asymptotic expansion of the form:
\begin{equation}\label{diagas}
\nu_k(p) \sim \kpi^n + \sum_{j=1}^{\infty}\, k^{n-j}\,a_{j}(p)
\end{equation}
For $v_h\in T_pZ$ a horizontal vector, we have the asymptotic
expansion:
\begin{equation}\label{pdp}
\langle \Psik(p)\,,\,d(\Psik)_p(v_h)\rangle \sim \sum_{j=0} c_j k^{n-1/2-j}.
\end{equation}
And given $u,v\in T_xX$, which lift to the horizontal vectors $u_h, v_h
\in T_pZ$, there is an asymptotic expansion of the form:
\begin{equation}\label{b1}
\langle d(\Psik)_p(v_h)\,,\,d(\Psik)_p(u_h)\rangle_{\calH_k}\, \sim \,
\frac{1}{2}\,
\frac{k^{n+1}}{(2\pi)^n}\,\bigl[ g(u,v) +\, i\, 
\omega (u,v) \bigr] +\sum_{j=0}^\infty\, c_j\,k^{n-j}.
\end{equation}
All three expansions are in the $C^\infty$ topology on the unit ball
bundle of $X$.
\end{theorem}

Applying this theorem to the formula (\ref{hxdef}) for $h_x$, we see 
that only the first term in formula for $h_x$ contributes to leading order.  The
pullback estimates in Theorem \ref{Main}, (\ref{dos}) and (\ref{tres}), then
follow immediately from (\ref{b1}).

\begin{proof}
The existence of the first asymptotic expansion (\ref{diagas}) 
is an immediate consequence of Lemma \ref{Tech}.  According to
Proposition \ref{szsym}, the symbol of $\Pi$ is a rank one projector
whose trace is 1,  so the leading term in the expansion is a constant
which depends only  on $k$ and the dimension.  The factor $(k/2\pi)^n$
may be computed  in a model case (the Bargmann kernel being the simplest
example).  The expansion is by integer steps because $\Pi$ is an
operator with definite (even) parity.
 
We turn next to (\ref{b1}).  Fix $p$ and $u_h,v_h\in T_pZ$ horizontal.
Extend $u_h,\,v_h$ to two  $S^{1}$-invariant
horizontal vector fields on $Z$, $U$ and $V$ respectively. 

\begin{lemma}\label{Lep}
Consider $U_1$, resp. $V_2$, as a differential operator on
$Z\times Z$ acting on the first, resp. second, variable.
Then
\begin{equation}\label{ugh}
\langle d(\Psik)_p(v_h)\,,\,d(\Psik)_p(u_h)\rangle_{\calH_k}\,=\, 
(U_1\circ V_2)\,(\Pi_{k}) (p,p).
\end{equation}
\end{lemma}
\begin{proof}
This follows from the reproducing property, (\ref{a1}).  
Let $p(s),\ q(t)$ be two
curves on $Z$ adapted to $u,\ v$ respectively.  Then
\begin{equation}\label{c4}
\langle\Psik(q(t))\,,\,\Psik(p(s))\rangle\, = \,\Pi_{k}(p(s), q(t)).
\end{equation}
Now just apply $\displaystyle{\frac{\partial^{2}}{\partial s\partial 
t}|_{(s,t)=(0,0)}}$.
\end{proof}

By the symbol calculus of Hermite distributions
(precisely applying Theorem 10.2 of \cite{BG}) one gets:
\begin{equation}\label{c1}
(U_1\circ V_2)(\Pi) \in J^{3/2}(Z\times Z, \Sigma)
\end{equation}
and its symbol at $(p, \alpha_{p}\,;\,p, -\alpha_{p})$
is known and will be described below. 
Notice furthermore that the right-hand side of (\ref{ugh}) is precisely
the $k$-th Fourier coefficient of $(U_1\circ V_2)(\Pi)$, in
the sense of Lemma \ref{Tech}.
Therefore (\ref{c1}) and Lemmas \ref{Tech} and \ref{Lep} imply 
the existence of the asymptotic expansion (\ref{b1}).
Moreover, $\Pi$ has definite parity and the application of pseudodifferential
operators preserves parity.  So the steps in the expansion are integral.

To compute the leading coefficient we must first 
describe the symbol of $(U_1\circ V_2)(\Pi)$.
Let $H_p\subset T_{p}Z$ be the horizontal subspace at $p$, and recall that
$u=U_{p}$, $v=V_{p}$.  Under the infinitesimal Heisenberg representation the 
vectors $u,v$ induce skew-hermitian operators, $\tilde{u},\ \tilde{v}$, 
on the smooth vectors of the metaplectic representation of $H_p$, $\calS_p$.
By Theorem 10.2 of \cite{BG} and the explicit form of the isomorphism
(\ref{iso}), the symbol of (\ref{c1}) at $p$ is
$|\tilde{u}(e)\rangle\otimes |\tilde{v}(e)\rangle$.  
We can compute this explicitly in terms of the K\"ahler form
on on the horizontal space $H_p$.

\begin{lemma}\label{Hrep}  
Let $(H,\, h)$ be a Hermitian vector space of dimension
$2n$ and $\calF$ a (Hilbert) representation space for the metaplectic
and Heisenberg representations of $H$.  Let $e\in\calF$ be a normalized
ground state of the harmonic oscillator.  Then:
\begin{equation}\label{c7}
\forall u,v\in H\qquad
\langle\tilde{u}(e)\,,\,\tilde{v}(e)\rangle_\calF\,=\,
\frac{1}{2}\,h(u,v).  
\end{equation}
\end{lemma}
\begin{proof}
It suffices to prove it for the standard model, $H=\bbC^{n}$ with
$h$ the standard hermitian form.  We can take:
$\calF = L^2(\bbR^n)$ and $ e\,=\,(\pi)^{-n/4}\,e^{-|x|^2/2}$.
If $u=a+ib$ where $a, b,\in\bbR^n$ then
$ \tilde{u}\,=\,ia\cdot x + b\cdot\nabla $
and so
\[
\tilde{u}(e^{-|x|^2/2})\,=\,
\bigl(i a\cdot x - b\cdot x \bigr)\,e^{-|x|^2/2}.
\]
A simple calculation yields (\ref{c7}).
\end{proof}

The expansion (\ref{b1}) now follows directly from Lemmas \ref{Tech} and
\ref{Hrep}.
The case of (\ref{pdp}) is quite similar.  Following Lemma \ref{Lep}, the
inner product in question may be written
$$
V_2\Pi \in J^1(Z\times Z,\Sigma).
$$
By Lemma \ref{Tech}, we get an asymptotic expansion $\sum_j c_jk^{n+1/2-j}$. 
A symbolic calculation similar to that of Lemma \ref{Hrep}
shows that the coefficient of the leading power, $k^{n+1/2}$, is zero.
\end{proof}

\subsection{$F_k$ is nearly pseudo-holomorphic}
Corollary \ref{Co} is proven as follows.
Let $J^0$ be the almost-complex structure on $\bbP\calH^*_k$.
At each point $p\in X$ we can decompose
$d(F_k)_p = \del (F_k)_p + \delb (F_k)_p$, where
$J^0 \circ \del (F_k)_p  = \del (F_k)_p \circ J$ and
$J^0 \circ \delb (F_k)_p  = - \delb (F_k)_p\circ J$.
Since $F_k$ is nearly symplectic and nearly an isometry, it is clear that
it should nearly intertwine the two almost-complex structures as well.
\begin{proposition}\label{Ac}
For $u \in T_pX$,
$$
\frac{1}{k} \norm{\del (F_k)_p u}^2 = \norm{u}^2 + O(k^{-1}), \qquad
\frac{1}{k} \norm{\delb (F_k)_p u}^2 = O(k^{-1}),
$$
and therefore
\begin{equation}\label{l1}
\frac{1}{k} \norm{ d(F_k)_p(J(u)) -  J^0(d(F_k)_p(u))}^2\,=\, O(k^{-1}).
\end{equation}
The estimates are uniform in $p$ and in $u$ on compact sets.
\end{proposition}

\begin{proof}
For notational simplicity introduce local coordinates and denote by $A$ 
the matrix of $d(F_k)_p$.  Then $\del (F_k)_p$
and $\delb (F_k)_p$ can be written as
$$
A^\pm = \frac12 \bigl[A\mp J^0AJ \bigr].
$$
So we need to compute
$$
\norm{A^\pm u}_{g^0}^2 = \frac14 \norm{Au}_{g^0}^2 + \frac14 \norm{J^0
A Ju}_{g^0}^2 \mp \frac12  g^0(Au\;,\; J^0 A J u).
$$
Applying Theorem \ref{Lapapa}, we find that
$$
\norm{A u}^2 = k\norm{u}_g^2 + O(1),
$$
and the same for $\norm{J^0 A Ju}^2$.
The remaining term is
\begin{equation*}
\begin{split}
g^0( A u\;,\; J^0 d(F_k)_p J u)
&=  \Omega_k(A u\;,\; A J u) \\
&=  k\omega(u, Ju) + O(1) \\
&=  k\norm{u}_g^2 + O(1). \\
\end{split}
\end{equation*}
\end{proof}

\subsection{Injectivity}

Theorem \ref{Lapapa} tell us that $F_k$ is an immersion
for large $k$. 
Consider the possibility that $F_k(p) = F_k(q)$.  This means that the
coherent state $\Psik(p)$ is a multiple of $\Psik(q)$.  Our first step is
to show that coherent states are concentrated at the base point for large
$k$.  So $q$ and $p$ would have to be asymptotically close together in order
for the coherent states to match.

In this subsection it will be convenient to denote by $d(\cdot,\cdot)$ not only
the distance function on $X$, but also its pullback to $Z$.

\begin{lemma}\label{leadterm}
For $p,q\in Z$,
$$
|\Pi_k(p,q)| = \nu_k(p) e^{-\frac{k}2 d(q,p)^2} + O(k^{n-1/2}),
$$
where the estimate is uniform in both $q$ and $p$.
\end{lemma}
\begin{proof}
One way to prove this is to observe that we can choose real phase functions 
$\theta_k(q,p)$ so that 
$$
G(q,p) = \sum_k \nu_k(p) e^{i\theta_k(q,p)} e^{-\frac{k}2 d(q,p)^2}
$$
is a distribution in $J^{1/2}(Z\times Z;\Sigma)$ with the same symbol
as $\Pi$.  (By taking Fourier transforms we can write $G$ in 
the form (\ref{ulocal}), and then read off its symbol and compare it to the harmonic
oscillator ground state.)
Thus $\Pi - G \in J^0(Z\times Z;\Sigma)$.  The uniform
estimate $|\Pi_k - G_k| = O(k^{n-1/2})$ is obtained from Lemma \ref{spest},
which in turn gives the estimate above.

Alternatively, one could directly argue as in Lemma \ref{spest}.  After obtaining
the leading term in $\Pi$ in (\ref{spresult}), one can simply compute $a_0(x,x,\eta)$ 
in local coordinates from the inverse of the symbol map.  
\end{proof}

\begin{corollary}\label{close}
There is a constant $C$ such that for any sequence $x_{m},y_{m} \in X$ such
that $F_{k_m}(x_{m}) = F_{k_m}(y_{m})$, $k_1<k_2<\dots$ we have
$$
d(x_{m},y_{m}) < Ck_m^{-3/4}.
$$
\end{corollary}
\begin{proof}
Choose $p_{m},q_{m}$ so that $\pi(p_m)=x_m$ and $\pi(q_m)=y_m$.  
By assumption
$\Psi_{k_m}(p_m) = \lambda_{m} \Psi_{k_m}(q_m)$.  By the reproducing property we have
$$
\nu_{k_m}(p_m) = |\lambda_m|^2\nu_{k_m}(q_m)
$$
so $|\lambda_m| = 1 + O(k_m^{-1})$. 
By Lemma \ref{leadterm} we have
$$
|\lambda_m|\> \nu_{k_m}(p_m) = |\Pi_{k_m}(p_m,q_m)| = \nu_m(p_m) 
e^{-\frac{k_m}2 d(p_m,q_m)^2} + O(k_m^{n-1/2}).
$$
This implies
$$
e^{-\frac{k_m}2 d(p_m,q_m)^2} \ge 1 - Ck_m^{-1/2}
$$
for some constant $C$.  By taking the logarithm and adjusting $C$,
we get 
$$
k_m d(p_m,q_m)^2 < Ck_m^{-1/2},
$$
and the result follows.
\end{proof}

To complete the proof, consider an arbitrary geodesic $\alpha:I \to X$,
parametrized by arclength.  We set $\gamma = F_k\circ \alpha:I\to
\bbP\calH_k$.  This will not of course be a geodesic, but we do have
uniform bound on the ``acceleration.''
\begin{lemma}\label{accel}
There is a constant $C$ such that for any $\alpha, k$, if $\gamma =
F_k\circ \alpha$ then
$$
\norm{\nabla_{\dot\gamma}\dot\gamma}_{g^0} < Ck.
$$
\end{lemma}
\begin{proof}
This is a matter of counting derivatives.  All terms in the expression for 
$\nabla_{\dot\gamma}\dot\gamma$ involve a total of 2
horizontal derivatives of $\Psik$.  
Let $L$ be an operator of degree $m$ given as
a product of horizontal vector fields on $Z$, which acts on $Z\times Z$ in the first 
variable.   Then $L\circ\Pi \in J^{(m+1)/2}(Z\times Z;\Sigma)$, and
by Lemma \ref{spest}, we have
$$
|(L\Pi_k)(p,q)| = O(k^{n+m/2}),
$$
uniformly in $q$ and $p$. The $k^n$ is cancelled by normalization, so to estimate
$\nabla_{\dot\gamma}\dot\gamma$
we just include a factor of $k^{1/2}$ for each derivative.
\end{proof}

Now we put these results together.
\begin{proposition}
The map $F_k$ is injective for all sufficiently large $k$.  
\end{proposition}
\begin{proof}
Suppose there is a sequence $x_m,y_m\in X$ such that $F_{k_m}(x_m) =
F_{k_m}(y_m)$ (with $k_1<k_2<\dots$).  By Lemma
\ref{close} we know $d(x_m,y_m)  = O(k_m^{-3/4})$.  For large $m$ let
$\alpha_{m}:I_{m} \to X$ be the geodesic segment connecting $x_m$ to $y_m$. 
Since $\norm{dF_{k_m}(\dot\alpha_{m})}_{g^0} = k_m^{1/2} + O(1)$, the curve
$\gamma_{m} = F_{k_m}\circ\alpha_{m}$ will be a closed ``lasso'' of 
length $O(k_m^{-1/4})$.

For each $m$, we'll choose normal
coordinates on $\bbP\calH_{k_m}$, centered at $F_{k_m}(x_m)$.
Let $|\cdot|$ denote the Euclidean norm in these coordinates.
If we draw a sphere centered at the origin and tangent to the lasso, 
at the point of tangency $\gamma(t_m)$ we have
$$
\frac{|\dot\gamma_m(t_m)|^2}{|\ddot\gamma_m(t_m)|} < R_m,
$$
where $R_m$ is the Euclidean radius of the sphere.

With normal coordinates the FS metric $g^0$ is the standard Euclidean metric to within
order $r^2$, and the length of 
the lasso is $O(k_m^{-1/4})$.  So in the neighborhoods of interest 
we can estimate the difference between the two metrics by $O(k_m^{-1/2})$.
Hence we can find a constant $C$ such that $R_m < Ck^{-1/4}$, and
we also have $|\dot\gamma(t_m)| = k_m^{1/2} + O(1)$.
The radius of curvature estimate becomes
$$
|\ddot\gamma_m(t_m)| > Ck_m^{5/4}.
$$

The acceleration vector is $\nabla_{\dot\gamma}\dot\gamma = \ddot\gamma + 
\Gamma\dot\gamma \dot\gamma$ (omiting the indices since they are of no concern here).
The Christoffel symbols can be bounded by a constant times $r$, 
so this can be estimated by 
$O(k_m^{-1/4})$, and each $\dot\gamma$ is $O(k_m^{1/2})$.
Thus $\norm{\nabla_{\dot\gamma}\dot\gamma(t_m)}_{g^0} = |\ddot\gamma(t_m)| + O(k^{3/4})$.
The estimate on $|\ddot\gamma_m|$ now gives
$$
\norm{\nabla_{\dot\gamma}\dot\gamma(t_m)}_{g^0} > Ck_m^{5/4},
$$
contradicting Lemma \ref{accel}.
\end{proof}

\section{Transverse hyperplane sections}

In this section we prove the following weak analogue
of a theorem of Donaldson:

\begin{theorem}\label{radonald}
Suppose $(X,\omega)$ is real-analytic (with $\omega$ giving an integral cohomology
class).  Then, for sufficiently large $k$, there exist 
hyperplanes $W$ transverse to $F_k$ such that $F_k^{-1}(W)$ is
symplectic off a codimension-one real analytic subset.
\end{theorem}

The proof is a consequence of the previous results and a
simple transversality argument.  We can choose the other
data, $J,L,\nabla$, to be real-analytic as well.  Then eigenfunctions
of the Laplacian $\Delta_k$ will be real-analytic so the coherent states
and the map $F_k$ will be real-analytic.   

We first need some preliminaries on dihedral angles.
Consider a complex $N$-dimensional Hermitian vector space, $\calV$,
let $W\subset \calV$ be a hyperplane and $V\subset\calV$ a complex
subspace of dimension $n$.    Let $w\in\calV$ be a unit
vector normal to $W$, and let $\pi_V : \calV\to V$ denote orthogonal
projection onto $V$.  The following is a measure of the dihedral
angle formed by $V$ and $W$:
\newcommand{\vt}{\vartheta}

\begin{definition}
We define
$\displaystyle{\vt (V,W) := \norm{\pi_V(w)}^2}$.
\end{definition}

Notice that $0\leq\vt(V,W)\leq 1$.  Also $\vt (V,W) = 0$ iff
$V\subset W$, and since $W$ is a hyperplane this occurs iff $V$ and
$W$ fail to intersect transversely.  And it should be clear that
$\vt (V,W) = 1$ iff $w\in V$, that is iff $W$ and $V$ intersect
orthogonally.

We can consider $\vt$ to be a function: 
\[
\vt:\ \calG_n \times \calG_{N-1}\to \bbR
\]
where $\calG_p$ denotes the Grassmannian of complex $p$-dimensional
subspaces of $\calV$.  As such $\vt$ is a real-analytic
function.

In our application, we will use $\vt$ to measure the angle between
a hyperplane $W\subset \bbP\calH_k^*$ and the complex subspace
$\del F_k(T_xX)$.  If $W$ is transverse to $F_k$ at $x$, 
then $V = T_x(F_k^{-1}(W)) \subset T_xX$ is a subspace of real codimension 2.
In fact, if $y = F_k(x)$ and $w \in T_{y}\bbP\calH_k^*$ is a unit orthogonal to 
$T_{y}W$, then $V$ is the kernel of the map $A:T_xX \to \bbC$, given by
$$
Au = \langle b, dF_k (u)\rangle.
$$
$A$ can be decomposed into its complex linear and antilinear parts: 
$A = A^+ + A^-$, where 
$$
A^+u = \langle b, \del F_k(u) \rangle, \qquad
A^-u = \langle b, \dbar F_k(u) \rangle.
$$
Recall from \cite{D} that $V$ is a symplectic subspace of $T_xX$ if
$|A^+|<|A^-|$ (the norm with respect to the hermitian structure on $T_xX$).

\begin{lemma}\label{apm}
\begin{equation*}
\begin{split}
|A^+|^2 &= k\>\vartheta(T_{y}W,\>\del F_k(T_xX)) + O(1), \\
|A^-|^2 &= O(1),  
\end{split}
\end{equation*}
where the estimates are uniform (in $x$).
\end{lemma}

\begin{proof}
This is a fairly direct consequence of Proposition \ref{Ac}.  The bound
on $A^-$ is immediate.  To estimate $A^+$, we use the fact that 
$\del F_k$ is $k$ times an isometry up to an error of $O(1)$.  
\end{proof}

Consider now $Y_k := F_k(X)\subset \bbP\calH_k^*$, with $k$ large.
Define the manifold
\begin{equation}\label{fl}
\calM := \{\, (y,W)\in Y\times \bbP\calH_k\;;\;y\in W\,\},
\end{equation}
where $W$ is to be thought of as a 
hyperplane in $\bbP\calH_k^*$. 
Consider the real-analytic projection
\begin{equation}\label{proj}
\begin{array}{ccc}
\calM & \to & \bbP\calH_k \\
(y, W) & \mapsto & W
\end{array}
\end{equation}
Applying Sard's theorem to the projection we obtain the existence
of many hyperplanes, $W$, intersecting $Y$ transversely:  their 
complement is a set of measure zero in the image of 
(\ref{proj}).

\begin{proof}[Proof of Theorem \ref{radonald}]
Given any $y\in Y$ there exists a hyperplane $W_0 \ni y$ such that 
$\displaystyle{  \vt(T_yW_0, \del F_k(T_x X)) = 1}$.  And we can
find $W$ arbitrarily close to $W_0$ such that $W$ is transverse to $Y$.
We can thus assume $W$ is transverse to $Y$ with 
angle $\vt(T_yW, \del F_k(T_x X))$ arbitrarily close to 1.  
If $k$ is sufficiently large then by Lemma \ref{apm} the restriction of $\omega$ to 
$F_k^{-1}(W) \subset X$ is non-degenerate at $x$.  

The determinant of the restriction of $\omega$ to
$F_k^{-1}(W)$ is a real-analytic function, which we know 
to be non-zero at a point.  Therefore
the set of zeroes is a real-analytic subset of codimension 1.
\end{proof}

\section{Toeplitz operators and dynamics}

Let $H:X\to \bbR$ be a smooth Hamiltonian.  To quantize
$H$ means to associate to it a sequence of self-adjoint
operators, $\{ T_k^H\}$, where
$T_k^H : \calH_k\to\calH_k$ for each $k$.  As already noted in
\cite{BU}, following ideas of Berezin \cite{B}, one way to define
$T_k^H$ is by the Toeplitz (or anti-Wick) prescription:
\begin{equation}\label{toe}
\forall\psi\in \calH_k\qquad T_k^H(\psi)\,=\,\Pi_k (H\psi),
\end{equation}
where $\Pi_k$ is the orthogonal projection onto $\calH_k$.
Since the projector $\Pi = \oplus \Pi_k$ defines a Toeplitz
structure in the sense of \cite{BG}, the assignment
$H\mapsto \{ T_k^H\}$ defines a deformation quantization,
\cite{BMS}, and the spectral estimates of \cite{BPU2} are
valid in the present setting as well.  The proofs of these
statements are identical to those in the K\"ahler case,
see {\em op.\ cit.}

Fixing a smooth Hamiltonian $H$ and suppressing the $H$-dependence
from the notation, we take $T_k := T_k^H$ to be the
quantum Hamiltonian corresponding to $H$.  The quantum 
dynamics are given by the sequence of 1-parameter subgroups,
\[
\tilde{U}_k(t)\,=\, e^{-ikt T_k} : \calH_k\to\calH_k,
\]
fundamental solution to the Schr\"odinger equation
$ik^{-1}\psi_t = T_k(\psi)$ (where $k = 1/\hbar$).
For each $t$ the unitary map $\tilde{U}_k(t)$ induces a
transformation of the projective space,
\[
U_k(t) : \bbP\calH_k \to \bbP\calH_k,
\]
which is holomorphic and an isometry.  On the other hand
one has the classical Hamiltonian flow of $H$, which we
will denote by $\phi_t : X\to X$.  The question arises:
to what extent are the embeddings $\Psik^o$ equivariant?
This is the issue we examine in this section.  

\newcommand{\Xik}{{}^k\Xi}
\newcommand{\tXik}{{}^k\tilde\Xi}

We begin by reviewing the Hamiltonian formulation of 
quantum mechanics.
Consider $\calH_k$ as a real symplectic vector space where the 
symplectic form is twice the imaginary part of the Hilbert inner product.
The natural $S^1$ action on $\calH_k$ is Hamiltonian with moment
map $\norm{\psi}^2$, and the projective space $\bbP\calH_k$ can be
thought of as a Marsden-Weinstein reduction of $\calH_k$ with
respect to this action.  

To be specific,
let $\tXik$ resp.\ $\Xik$ denote the infinitesimal generator of 
the one-parameter group $\tilde{U}(t)_k$ resp.\ $U_k(t)$.  
$\tXik$ is the linear vector field on $\calH_k$
\begin{equation}\label{campo}
\forall \psi\in \calH_k\qquad
\tXik_\psi = -ik\, T_k \psi,
\end{equation}
which is the Hamiltonian vector field of the function
\[
\tilde{Q}_k(\psi) = k\,\langle T_k \psi\, ,\, \psi\rangle.
\] 
The Hamiltonian is invariant under the $S^1$ action
and $U_k$ is the reduction of $\tilde{U}_k$ to $\bbP\calH_k$.  
In other words, the vector field $\Xik$ is the Hamiltonian vector field of 
the function \begin{equation}\label{wi}
Q_k([\psi]) = k\,\frac{\langle T_k \psi\,,\,\psi\rangle}
{\langle\psi\,,\,\psi\rangle}
\end{equation}
(where $[\psi]$ denotes the complex line through $\psi\not= 0$) with
respect to the  Fubini-Study symplectic form.
This is known, see \cite{AS} for example and references therein.
In addition to the symplectic structure, in quantum mechanics
one also has a Riemannian metric on the phase space, namely the
Fubini-Study metric on $\bbP\calH_k$.  This has a nice physical interpretation:
\begin{lemma}\label{Qlength}
The length of $\Xik$ at the quantum state $[\psi]$ is 
\[
|\Xik_{[\psi]}|^2 =
2k^2\,\frac{\langle T^2_k \psi\,,\,\psi\rangle}{\langle\psi\,,\,\psi\rangle} - 
2k^2\,\left(
\frac{\langle T_k \psi\,,\,\psi\rangle}{\langle\psi\,,\,\psi\rangle}\right)^2,
\]
That is, up to a factor of $2k^2 = 2/\hbar^2$,
it is the mean-square deviation in the observation of the
energy when the quantum system is in the state $[\psi ]$.
\end{lemma}

\newcommand{\sbot}{\bot_s}

We can now state some partial results on near-equivariance of $\Psik^o$.

Consider the symplectic submanifold $Y_k := \Psik^o (X)$.
(Although $\Psik^o$ is not exactly symplectic its image is a symplectic 
submanifold of $\bbP\calH_k$.)  
If $\sbot$ denotes the symplectic orthogonal, one has:
\begin{equation}\label{sbot}
\forall y\in Y_k \qquad T_y(\bbP\calH_k) = T_yY_k \oplus (T_yY_k)^{\sbot}.
\end{equation}
We will use the notation: $v = v^{\|} + v^{\sbot}$ for the 
components of a vector $v$ under this decomposition.  In particular,
the vector field along $Y_k$, $\Xik|_{Y_k}$, splits under this 
decomposition
\[
\Xik|_{Y_k}\,=\, \Xik|_{Y_k}^{\|} + 
\Xik|_{Y_k}^{\sbot}.
\]
We will show that the tangential component is
approximately equal to $d\Psik^o(\xi_H)$:

\begin{theorem}\label{Equidetail}
One has, uniformly,
\begin{equation}\label{uxmal}
d\Psik^o(\xi_H) \,=\, \Xik |_{Y_k}^{\|} + O(k^{-1/2}).
\end{equation}
In particular $ \Xik |_{Y_k}^{\|} = O(k^{1/2})$.
Moreover, assuming that $\Xik^{\sbot}\not=\vec{0}$,
let $\theta_k (y)$ denote the angle between $\Xik^{\|}$ and 
$\Xik^{\sbot}$ at $y\in Y_k$.  Then
\begin{equation}\label{langle}
|\cos (\theta_k(y))| = O(k^{-1/2})
\end{equation}
uniformly on $y$.  If $X$ is K\"ahler then $\cos (\theta_k)$ is 
identically zero.
\end{theorem}

Before giving the proof we'll go over a few facts.
Consider the inclusion map
\[
\iota : Y_k\hookrightarrow \bbP\calH_k .
\]
The following is immediate:
{\em 
The Hamiltonian vector field of the restriction $\iota^* Q_k$ of $Q_k$
to $Y_k$ is the tangential component of $\Xik$ with respect to the
above decomposition:}
\begin{equation}\label{obv}
\Xi_{\iota^* Q_k} = \Xik |_{Y_k}^{\|}
\end{equation}
\newcommand{\wick}{\sigma_{T_k}}
The function on $X$
\begin{equation}
\wick (x) := \frac{1}{k}\,\Psik^{o*}Q_k (x) = 
\frac{\langle T_k(\Psi_k (p))\,,\,\Psi_k (p)\rangle}
{\langle\Psi_k (p)\,,\,\Psi_k (p)\rangle}
\end{equation}
(where $p\in Z$ sits above $x$) is called the {\em Wick} or
{\em covariant}, \cite{B}, symbol of the operator $T_k$.
Recall that $T_k$ is the level $k$ Toeplitz operator associated
to a smooth Hamiltonian $H : X\to\bbR$.  We will refer to $H$ as
the Toeplitz symbol of $T_k$.

\begin{lemma}\label{Covest} 
As $k\to\infty$ the Wick symbol has an asymptotic expansion
of the form
\[
\wick \sim H + \sum_{j=1}^\infty\; \sigma_j\ k^{-j}.
\]
The asymptotics are in the $C^\infty$ topology.
In other words, to leading order the Wick and Toeplitz symbols agree.
\end{lemma}
\begin{proof}
This is another calculation in the Hermite calculus, using exactly
the same technique used in the calculation of the asymptotics of the
function $\nu_k$ in Theorem \ref{Lapapa}.
\end{proof}

\begin{proof}[Proof of Theorem \ref{Equidetail}]
By (\ref{obv}) the pull-back, $\xi_k := \Psik^{o*} (\Xik |_{Y_k}^{\|})$ 
is the Hamilton vector field of $k\wick$ with respect to the symplectic form 
$\omega_k:= \Psik^{o*}\Omega_k$:
\[
k\,d(\wick )\, =\, \omega_k \rfloor \xi_k,
\]
and by Lemma \ref{Covest}
this is $O(k)$ with respect to the fixed Riemannian metric $g$.
On the other hand, by Theorem \ref{Main}, $\omega_k = k\omega + O(1)$.
Therefore $\xi_k = O(1)$. 
But by Lemma \ref{Covest} $d(\wick ) = dH + O(1/k)$.
Therefore
\[
k\,dH + O(1) = k\omega\rfloor\xi_k + O(1)\rfloor\xi_k
\]
and since the last term is $O(1)$ we can conclude
$dH = \omega\rfloor\xi_k + O(1/k)$.  Therefore
\begin{equation}\label{hquasi}
\xi_H = \xi_k + O(1/k),
\end{equation}
where this estimate is with respect to the fixed Riemannian metric
$g$ on $X$.  Applying $d\Psik^o$ and using Theorem \ref{Main} again gives
(\ref{uxmal}).

To prove (\ref{langle}) we will estimate the Euclidean inner product
\[
g_k^0(\Xik^{\|}\,,\,\Xik^{\sbot})\,=\,-\Omega_k(J^0\,\Xik^{\|}\,,\,\Xik^{\sbot}),
\]
where $J^0$ is the complex structure on the real tangent bundle of
$\bbP\calH_k$.
The idea is that $\Psik^o$ is nearly antiholomorphic, so 
$J^0(\Xik^{\|})$ is nearly tangent to $Y_k$.  More precisely,
by (\ref{l1}) we have 
\[
d\Psik^o(\nabla H) = d\Psik^o(J\xi_H) = -J^0\,d\Psik^o(\xi_H) + O(1)
\]
Coupled with (\ref{uxmal}), we get:
\[
J^0(\Xik^{\|}) = -d\Psik^o(\nabla H) + \epsilon_k
\]
where $|\epsilon_k| = O(1)$. Since $d\Psik^o(\nabla H)$ is tangent
to $Y_k$ and $\Xik^{\sbot}$ is in the symplectic orthogonal,
\[
g_k^0(\Xik^{\|}\,,\,\Xik^{\sbot})\,=\ -\Omega_k 
(\epsilon_k\,,\,\Xik^{\sbot})\,=\,g_k^0(\epsilon_k\,,\,J^0(\Xik^{\sbot})).
\]
By the definition of $\theta_k$, Schwartz' inequality and the fact that
$J^0$ preserves the metric, we get
\[
|\Xik^{\|}|\,|\Xik^{\sbot}|\,|\cos (\theta_k)| \leq 
|\Xik^{\sbot}|\,|\epsilon_k|.
\]
Recalling that $|\epsilon_k|= O(1)$ and that $|\Xik^{\|}| = O(k^{1/2})$
yields the result.
\end{proof}

\bigskip\noindent
{\bf Final remarks.}
In several special cases (e.\ g.\ $X=\bbC\bbP^1$)
we can prove, by direct calculation, that for all $H$
\begin{equation}\label{uxmal2}
\Xik|_{Y_k}^{\sbot}\,=\,O(1), 
\end{equation}
from which the estimate
\begin{equation}\label{bett}
 d(\Psik^o)_x (\xi_H) \,=\, \Xik_{\Psik^o(x)} + O(1)
\end{equation}
follows.   The estimate (\ref{bett}) is an infinitesimal 
version of the near-equivariance of the $\Psik^o$.
(Notice that by Theorem \ref{Main} $d\Psik^o(\xi_H) = O(k^{1/2})$.)
By Lemma \ref{Qlength},
the issue of whether (\ref{uxmal2}) holds is equivalent to
the issue of whether
\begin{equation}\label{hope}
\frac{\langle T^2_k(\Psi(p))\,,\,\Psi(p)\rangle}{\langle\Psi(p)\,,\,\Psi(p)\rangle} -
\left( \frac{\langle T_k(\Psi(p))\,,\,\Psi(p)\rangle}{\langle\Psi(p)\,,\,\Psi(p)\rangle}\right)^2 =
|\left(\xi_H\right)_x|^2 \,k^{-1} + O(k^{-2})
\end{equation}
holds.  (Here $p\in Z$ is a point above $x\in X$.)
We find this to be an interesting semi-classical question: 
it relates the mean-square deviation of the energy, when the
system is in a coherent state, to the Riemannian length of the
classical Hamiltonian vector field.
One can easily see that for a given $(X,\omega, J)$ the leading
asymptotics of the left-hand side of (\ref{hope})
are unchanged if one modifies $T=\Pi H \Pi$ by a
pseudodifferential operator of order $(-1)$.
This supports the conjecture that 
(\ref{hope}), and therefore (\ref{bett}), hold in general.
We hope to return to this issue in the near future.
(It is easy to check that (\ref{hope}) holds for $X=\bbC^n$ with the
Euclidean metric.)

\appendix

\section{The spectral projector}

We will prove here Theorem \ref{Bg} and Proposition \ref{szsym}, that is,
that the projector $\Pi : L^2(Z)\to\calH$ is an Hermite FIO of the same
form as the Szego projector. 
Our starting point is the following result from \cite{GU}:

\begin{theorem} \cite{GU}
There exists a self-adjoint second-order pseudodifferential operator,
$Q$ on $Z$, commuting with $S^1$ such that:
\begin{enumerate}
\item The orthogonal projector, $S:L^2(Z)\to \calW$ onto
the $L^2$ closure of the kernel of $Q$, is an Hermite FIO 
in the class $S\in J^{1/2}(Z\times Z, \Sigma)$
and symbol as claimed in Proposition \ref{szsym} for $\Pi$.
\item The non-zero spectrum of $Q$ is positive, and it drifts in the 
same sense as that of $\calB$.
\item One has: $\calB\,=\,Q + R$, where $R$ is a classical zeroth-order
$\Psi$DO commuting with $D_\theta$.
\end{enumerate}
\end{theorem}

To prove Theorem \ref{Bg}, we will show that $\Pi - S$ is a smoothing
operator.  We will use the notation $A\sim B$ if $A-B$ is smoothing.

\begin{lemma}\label{rlemma}
Let $R$ be a classical zeroth-order $\Psi$DO, commuting with $D_\theta$.
There exists a another pseudodifferential operator $T$, commuting with 
$D_\theta$, such that
$$
RS \sim TS \qquad\text{and}\qquad [T,S] \sim 0.
$$
\end{lemma}
\begin{proof}
Given any zeroth-order $R$ there exists, by Proposition 2.13 of \cite{BG}, 
a zeroth-order operator $T_0$
which commutes with $S$ and $D_\theta$, so that $ST_0S = SRS$.  This means
that $\sigma(R-T_0)$ must vanish on 
$\Sigma$.  Applying Theorem 10.2 of \cite{BG}, 
we see that $(R-T_0)S$ is an Hermite FIO of
order $-1/2$.  Corollary 5.6 of the Appendix of \cite{BG} then implies that
there exists a $\Psi$DO of order $-1/2$, $R_{-1/2}$, such that 
$$
(R-T_0)S = R_{-1/2}S.
$$

At the next stage, we start with $R_{-1/2}$, construct $T_{-1/2}$ which
commutes with $S$ so that $(T_{-1/2}-R_{-1/2})S$ is order $-1$, and so on.

We sum the resulting sequence $T_j$
asymptotically to get $T$ which commutes with $S$ up to smoothing,
and for which $(R-T)S \sim 0$.
\end{proof}

The following Proposition clearly implies Theorem \ref{Bg}.
\begin{proposition} $\Pi - S$ is a smoothing operator:
$$
\Pi \sim S.
$$
\end{proposition}
\begin{proof}
We will use the identity (valid for any two projections)
\[
\Pi - S\,=\,\Pi S^\bot - \Pi^\bot S
\]
and prove that each term on the right-hand side is smoothing.

Let's prove first that $\Pi^\bot S$ is smoothing.  For some $N$,
decompose 
$L^2(Z) =  L^2(Z)_{\ge N} \oplus L^2(Z)_{<N}$,
where $L^2(Z)_N = \bigoplus_{k=N}^\infty L^2(Z)_k$. 
Since the wave front set of $S$ is $\Sigma$ we can deduce
immediately that $S|_{L^2(Z)_{<N}}$ is smoothing for any $N$.

So we turn our attention to $L^2(Z)_{\ge N}$.  By the hypoellipticity
Theorem of Boutet de Monvel if $\epsilon >0$ is small
$\calB+\epsilon D_\theta$ is hypoelliptic with loss of one derivative.
By the definition of $\calH$ and the drift of the spectrum of $\calB$,
if $N$ is sufficiently large there is a constant $C$ such that for any
$u_k\in\calH^\bot\cap L^2(Z)_k$ with $k>N$,
$$
\norm{\calB u_k} \ge Ck\norm{u_k}.
$$
But $k\norm{u_k} = \norm{D_\theta u_k}$, so $D_\theta$ is dominated by $C\calB$
on $\calH^\bot\cap  L^2(Z)_{\ge N}$: 
\[
\exists C\ \forall u\in\calH^\bot\cap L^2(Z)_{\ge N}\qquad
\norm{ D_\theta u} \leq C\norm{\calB u}.
\]
So (recall that both $\calB$ and $Q$ commute with $D_\theta$)
\[
\forall u\in\calH^\bot\cap L^2(Z)_{\ge N}\ \forall m\geq 0\qquad
\norm{ (\calB+\epsilon D_\theta)^m u} \leq C_{m,\epsilon}\,\norm{\calB^m u}.
\]
Therefore, if we can prove that 
$\calB^m\Pi^\bot S$ is bounded on $L^2(Z)_{\ge N}$ for each $m\ge 0$, 
then the hypoellipticity of
$\calB+\epsilon D_\theta$ will imply that $\Pi^\bot S$ is smoothing on
$L^2(Z)_{\ge N}$.  

According to Lemma \ref{rlemma}, $\calB = Q+R \sim Q + TS + RS^\bot$ with
$[T,S]=0$.  Therefore (using that $\calB S = \calB (\Pi + \Pi^\bot ) S$)
$$
\calB \Pi^\bot S = \calB S - \calB \Pi S \sim TS - \calB \Pi S.
$$
Since $T$ commutes with $S$ up to smoothing, after repeated
applications of $\calB$ we obtain
$$
\calB^m \Pi^\bot S \sim T^m S - \calB^m \Pi S.
$$
Since $T$ is order zero, $T^mS$ is a bounded operator.  $\calB^m \Pi$ is
also bounded, by the definition of $\calH$.
Thus $\calB^m \Pi^\bot S$ is bounded (on all of $L^2(Z)$ in fact),
and we have $\Pi^\bot S \sim 0$ by hypoellipticity.

To show that $\Pi S^\bot$ is smoothing, we'll show that
$(\Pi S^\bot)^* =S^\bot \Pi$ is smoothing.
Since the part of the spectrum of $Q$ corresponding to the image of $S^\bot$
drifts just like that of $\calB$,
we know that $D_\theta$ is dominated by $Q$ on the image of $S^\bot$
intersection $L^2(Z)_{\ge N}$.
Since $\calB = Q+R$ with $R$ bounded, we can make the estimate
that $\exists C>0$ such that for sufficiently large $N$ 
\[
\forall u\in L^2(Z)_{\ge N}\qquad \norm{\calB S^\bot u } >
C \norm {D_\theta S^\bot u}.
\]
Therefore, by the same reasoning as above,
it suffices to show that $\calB^m S^\bot \Pi$ is bounded for all $m\ge 0$.
But we have:
$$
\calB S^\bot \Pi = \calB(1-S)\Pi = \calB\Pi-(Q+R)S\Pi = \calB\Pi-RS\Pi.
$$
By the first part of the proof, $\Pi^\bot S \sim 0$, or in other words
$S \sim \Pi S$.   Since $S$ and $\Pi$ are self-adjoint we then have $S\sim S 
\Pi$
also.  Hence $RS\Pi \sim RS \sim TS$ and 
$$
\calB S^\bot \Pi \sim \calB\Pi - TS.
$$
Iterating this computation, we obtain
$$
\calB^m S^\bot \Pi \sim \calB^m\Pi-T^m S,
$$
which is bounded.  So $S^\bot \Pi \sim 0$ and the theorem is proved.
\end{proof}

\section{The symbolic calculation}
\newcommand{\hffm}{{\textstyle \bigwedge^{1/2}}}
\newcommand{\nchffm}{{\textstyle \overline{\bigwedge}^{-1/2}}}

We begin by briefly reviewing the symplectic linear algebra that underlies
the symbolic calculation in Lemma \ref{Tech}.
Let $V$ and $W$ be symplectic vector spaces, $\Gamma\subset V\times W^{-}$
a Lagrangian subspace and $\Sigma\subset W$ an isotropic subspace.  We
think of $\Gamma$ as a canonical relation from $W$ to $V$; $\Gamma\circ\Sigma$
is an isotropic subspace of $V$.  (For the composition $Pu$ in Lemma \ref{Tech},
$\Gamma$ and $\Sigma$ are the 
linearization of the canonical relations of $P$ and $u$, respectively.)
Define the symplectic vector space
\begin{equation}\label{appa}
\calN\,:=\,\Sigma^{\bot}/\Sigma.
\end{equation}
A symplectic spinor on $\Sigma$ is an element of $\calS(\calN)\otimes \hffm (\Sigma)$.
Since $\Gamma$ is Lagrangian, a symbol on $\Gamma$ is just a half-form in 
$\hffm (\Gamma)$.

The general symbol map for composition of a Lagrangian distribution with
a Hermite distribution is given in \S6 and \S7 of \cite{BG}.
Here we will present a simpler composition formula which holds under
two key assumptions:
\begin{equation}\label{app1}
\Gamma\circ\Sigma\ \mbox{is Lagrangian}
\end{equation}
and
\begin{equation}\label{app2}
\{\,w\in\Sigma\,;\,(0,w)\in\Gamma\,\}\,\subset W\,=\,\{ 0\} .
\end{equation}
(Both assumptions will hold for the composition of $P$ with $u$.)
In general, the result of composition will be a symbol on $\Gamma\circ\Sigma$.  
Since we assume that this space is
Lagrangian, the symbol will just a half-form in $\hffm (\Gamma\circ\Sigma)$.
Thus the symbol composition formula we are looking for will be a map
\begin{equation}\label{app3}
\calS(\calN)\otimes \hffm (\Sigma) \otimes \hffm (\Gamma)
\to \hffm ( \Gamma\circ\Sigma ).
\end{equation} 

Of central importance in the symbol composition is the subspace
\[
U_1\,:=\,\{\,w\in\Sigma^{\bot}\,;\,(0,w)\in\Gamma\,\}\,\subset W.
\]
We will denote by $U$ the image of $U_1$ in $\calN$.
From (\ref{app2}) it follows that the projection $U_1\to U$ is an
isomorphism.

Consider
\[
\begin{array}{rcc}
\rho : \Gamma \oplus \Sigma  & \to &  W \\
((v,w),w_1) & \to & w-w_1
\end{array}
\]
One can show that the image of this map is exactly $U_1^{\bot}$, and
under the assumption (\ref{app2}) its kernel is isomorphic to
$\Gamma\circ\Sigma$.  In other words one has an exact sequence:
\begin{equation}\label{exactseq}
0 \to \Gamma\circ\Sigma \cong \ker (\rho)  \to \Gamma \oplus \Sigma
\to U_1^{\bot}\to  0.
\end{equation}

Using this exact sequence and (\ref{app1}), we can see that $U\subset \calN$ is a 
Lagrangian subspace, as follows.
First note that $U$ will always be isotropic, because $\Gamma$ is isotropic.  
By (\ref{exactseq}), 
$$
\dim U_1^\bot = \frac12(\dim V + \dim W) + \dim \Sigma - \dim \Gamma\circ\Sigma.
$$
But $ \dim\Gamma\circ\Sigma = \frac12 \dim V$ by (\ref{app1}), so
$$
\dim U_1^\bot = \frac12 \dim W + \dim \Sigma.
$$
And this means 
$$
\dim U = \dim U_1 =  \frac12 \dim W - \dim \Sigma = \frac12 \dim\calN,
$$
so $U$ is Lagrangian.

We are now prepared to describe the two ingredients which yield the symbol
map (\ref{app3}).

\begin{lemma}\label{hfmap}  Under the assumption (\ref{app2})
there exists a canonical isomorphism
\[
\hffm (\Sigma)\otimes \hffm (\Gamma)\;\cong\;
\hffm (U) \otimes \hffm ( \Gamma\circ\Sigma ).
\]
\end{lemma}
\begin{proof}
When applied to the exact sequence (\ref{exactseq}), the functor $\hffm$
gives
\[
\hffm (\Sigma)\otimes \hffm (\Gamma)\;\cong\;
\hffm (U_1^\bot) \otimes \hffm ( \Gamma\circ\Sigma ).
\]
We also have an isomorphism $U_1^\bot \cong U \oplus (U_1^\bot/U_1)$.
The space $U_1^\bot/U_1$ is symplectic and carries a natural half-form.
The final result is obtained by dividing out by this form.  
\end{proof}

The second ingredient is the Kostant pairing:

\begin{lemma}\label{kostant}(Kostant)  If $\calN$ is a symplectic vector space 
and
$U\subset \calN$ a Lagrangian subspace there is a natural pairing
\begin{equation}\label{app4}
\calS(\calN) \otimes \hffm (U) \to \bbC
\end{equation}
\end{lemma}
\begin{proof}
$U$ is a maximal abelian subalgebra of the Heisenberg group of $\calN$,
and under the Heisenberg representation the intersection of the kernels 
of the operators in $U$ form a 1-dimensional subspace of
$H_{-\infty}$, the dual space to $\calS(\calN)$ (space of 
tempered distributions).  Kostant proved that this subspace is
naturally isomorphic with $\hffm (U)$.
\end{proof}

To obtain the symbol map (\ref{app3}), apply first the map from Lemma 
\ref{hfmap} and then the pairing from Lemma \ref{kostant}.

\medskip
We now finish the symbolic calculation in Lemma \ref{Tech}.  
The operator $P$ from (\ref{pfio})
is a Fourier integral operator associated to the following canonical
relation, $\Gamma\subset T^*(S^1\times Z\times Z)$:
\[
\Gamma\,=\, \{\, (\theta, <\eta,\dtheta>\,;\,p\cdot e^{i\theta}, \eta\,;\,
p, \xi)\,:\theta\in S^1,\; \eta,\xi\in T_pZ\}.
\]
To define the symbol of $P$, we must choose a trivialization of the half-form 
bundle of $\Gamma$
consistent with that chosen for $\Sigma$ in \S3 
(because after the half-forms are divided
out we need to have $Pu(\theta) = u(p\cdot e^{i\theta}, p)$).
The symbol of $P$ is then the combination of the half-form on $T_p^*P$ derived 
from
the metric and the half-form on $S^1$ determined by the standard metric. 
We also equip $\Gamma \circ \Sigma = T_0^*S^1$ with the natural half-form
given by this identification.

Let $Y = {(p;r\alpha_p); p\in Z, r>0} \in T^*P$, a symplectic submanifold.
Then from (\ref{1d}) the isotropic $\Sigma$ is the diagonal in 
$Y\times Y^-$, where the minus indicates the symplectic form has been reversed.  
Thus the symplectic normal $\calN = Y^\bot \times (Y^\bot)^-$, 
and $U$ is the diagonal in this space.   Since both $U$ and
$\Sigma$ carry symplectic structures, they possess natural half-forms
as well.  In \S3, we defined the symbol of $u \in J^m(Z\times Z; \Sigma)$ 
using we are using the natural half-form on $\Sigma$. 

The following shows that with these natural trivializations the half-form part 
of the symbol map drops out.
\begin{lemma}
In this case, the map of Lemma \ref{hfmap} takes the natural half-forms on the
right-hand side to the natural half-forms on the left.
\end{lemma}
\begin{proof}
We linearize the problem at the points 
$$
(0,1;p,\alpha_p; p,-\alpha_p) \in \Gamma, \qquad
(0,1;p,\alpha_p; p,-\alpha_p)\in\Sigma,
$$ 
and corresponding points in the other
spaces.  To avoid cluttered notation, for the remainder of this proof we'll use
the same letters $\Gamma,\Sigma, \calN, Y$, etc. to denote the 
{\it tangent spaces} to these spaces at the appropriate points.

Let $L = \{\xi\in T_{\alpha_p}(T_p^*Z):\; \xi(\dtheta) = 0\}$.
As a relation, $L$ defines a map from $Y^\bot$ to $Y$.
That is, given $w\in Y^\bot$, there is a unique $v\in Y$ such that
$w+v\in L$.  We'll denote this map by $w\mapsto L(w)$.

Under the projection $T(T^*Z) \to TZ$, the space $Y^\bot$ projects 
to the horizontal subspace of $T_pZ$, which we'll denote by $H$.  
This map is a symplectomorphism.  The space $L$ may also be identified with
$H$, under the map $w\mapsto w + L(w)$.  Define $V$ to be the span of 
$\{\dtheta;\alpha_p\} \in T(T^*Z)$, a 2-dimensional symplectic vector
space.   Then using the map $L$, we obtain a symplectomorphism
$Y \cong H \oplus V$.  Thus $\Sigma$ is also identified with $H\oplus V$.

Let $\Lambda = \Gamma\circ\Sigma = T_1(T_0^*S^1)$. 
Then $\Gamma$ is identified with $\Lambda \oplus H \oplus H \oplus V$,
using 
$$
(\sigma,w,w',u) \mapsto (0,\sigma; w+L(w) + u + \sigma \alpha_p; 
w'+L(w') - \sigma \alpha_p) \in \Gamma.
$$
Under this map, the natural half-form on $\Gamma$ is the combination of the 
natural
half-forms on $\Lambda, H$, and $V$.
 
We will split the exact sequence (\ref{exactseq}) into 3 parts, according to 
these decompositions. The first, and simplest, is just
$$
0\to\Lambda\to\Lambda\to0\to0.
$$
The second $\Lambda$ is identified with a subspace of $\Gamma\oplus\Sigma$
under $\sigma\mapsto (0,\sigma; \sigma\alpha_p; -\sigma\alpha_p)\oplus
(\sigma\alpha_p; -\sigma\alpha_p)$.  The corresponding half-form map is just 
the identity.

The second part consists of the horizontal spaces:
$$
0\to 0\to H\oplus H \oplus H \to  H\oplus H \oplus H \to 0.
$$
Here the first $H\oplus H \oplus H$ is identified with a subspace of $\Gamma\oplus\Sigma$
by
$$
(u,v,w) \mapsto (0,0;u+L(u); v+L(v))\oplus(L(w);L(w)) \in \Gamma\oplus\Sigma,
$$
and the second we make the identification
$$
(u,v,w) \mapsto (L(u);L(u))\oplus (v+L(v)+w;w) \in U\oplus (U_1^\bot/U_1).
$$
The resulting map is the symplectomorphism $(u,v,w)\mapsto (v-w,u-v,v)$,
so the half-form map preserves the natural half-forms.

The final part of the sequence is
$$
0\to0\to V\oplus V\to V\oplus V\to0.
$$
The first identification is
$$
(u,v) \mapsto (0,0;u;0)\oplus(v;v) \in \Gamma\oplus\Sigma,
$$
and the second is the inclusion $V\oplus V \subset U_1^\bot$.
The resulting map is $(u,v)\mapsto (u-v,-v)$, again a symplectomorphism.
\end{proof}

\medskip
To complete the proof of Lemma \ref{Tech}, 
consider the Kostant pairing for the case where $\calN = 
H_p \times (H_p)^-$, and $U$ is the diagonal. 
The distribution associated to the natural half-form on $U$ is just a delta-function
on the diagonal, so the Kostant pairing gives integration on the diagonal. 
If we write the element of $\calS(\calN)$ as a map $\sigma:\calS(H_p)\to \calS(H_p)$, 
we just get the trace of $\sigma$, as claimed in Lemma \ref{Tech}.

\end{document}